\documentclass[12pt]{article}
\usepackage{amsmath,amssymb,amsthm}
\usepackage{graphicx,psfrag,float,subcaption,caption}
\usepackage{enumerate,pifont,mathrsfs}
\usepackage[normalem]{ulem}
\usepackage[scale=0.8,a4paper]{geometry}
\usepackage[colorlinks,linkcolor=blue,citecolor=blue]{hyperref}
\usepackage{tikz}

\newtheorem{theorem}{Theorem}[section]

\newtheorem{lemma}[theorem]{Lemma}
\newtheorem{corollary}[theorem]{Corollary}

\newtheorem{prop}[theorem]{Proposition}

\newtheorem*{observation}{Observation}

\title{Claw-free minimal factor-critical and bicritical graphs}
\author{
{\small\bf Huaibin Zhou, Mingyang Guo and Xiumei Wang\thanks{Corresponding author. E-mail: wangxiumei@zzu.edu.cn}}\\
{\small School of Mathematics and Statistics, Zhengzhou University,}
{\small Zhengzhou, Henan 450001, China}
}
\date{}

\begin{document}
\maketitle

\begin{abstract}
A graph $G$ is called $k$-factor-critical if $G-S$ has a perfect matching for every subset $S\subseteq V(G)$ with $|S|=k$, and is minimal if $G-e$ is not $k$-factor-critical for every edge $e$ of $G$.
The cases $k=1$ and $k=2$ correspond to factor-critical and bicritical graphs, respectively.
In this paper, we obtain complete structural characterizations of claw-free minimal factor-critical graphs and claw-free minimal bicritical graphs.
The former is characterized by ear decompositions.
For the latter, we first characterize the $3$-connected case using the theory of line graphs and then extend the characterization to the general case by a suitable $2$-sum operation.


\end{abstract}

\section{Introduction}
All graphs considered in this paper are finite, connected, undirected, and simple unless otherwise stated.
We use the standard terminology and notation of \cite{BM08,LovaszPlummer}.
For a graph $G$, let $V(G)$ and $E(G)$ denote its vertex set and edge set, respectively.
The order of $G$ is $|V(G)|$.
For $S\subseteq V(G)$, let $G[S]$ denote the subgraph of $G$ induced by $S$.
For a vertex $v\in V(G)$, let $N_G(v)$ denote the set of neighbours of $v$ in $G$, and the degree of $v$, denoted by $d_G(v)$, is the number of edges of $G$ incident with $v$.
A \emph{matching} of $G$ is a set of pairwise nonadjacent edges.
A matching is \emph{perfect} if it covers every vertex of $G$.
For an integer $k$ with $1\leq k\leq |V(G)|-1$, a graph $G$ is \emph{$k$-factor-critical} if $G-S$ has a perfect matching for every subset $S\subseteq V(G)$ with $|S|=k$.
In particular, $1$-factor-critical and $2$-factor-critical graphs are called \emph{factor-critical} and \emph{bicritical}, respectively.

Factor-critical and bicritical graphs play important roles in matching theory.
Factor-critical graphs arise naturally in the Gallai-Edmonds decomposition of maximum matchings \cite{Edmonds1965,Gallai1964,LovaszPlummer} and admit a structural characterization in terms of odd ear decompositions due to Lov\'asz \cite{Lovasz1972}.
Bicritical graphs are closely related to bricks: a graph is a brick if and only if it is $3$-connected and bicritical \cite{EdmondsLovaszPulleyblank}.
Bricks are fundamental building blocks in the tight cut decomposition of matching covered graphs \cite{Lovasz1987}.

A graph is \emph{claw-free} if it contains no induced subgraph isomorphic to $K_{1,3}$.
Claw-free graphs are an important class in structural graph theory and matching theory, including several well-studied graph families such as line graphs and circular interval graphs \cite{Faudree1997}.
Chudnovsky and Seymour \cite{ChudnovskySeymour2008} provided a structural characterization of claw-free graphs.
Bruhn and Stein \cite{BruhnStein2012} characterized claw-free $t$-perfect graphs, that is, claw-free graphs whose stable set polytopes are described by nonnegativity, edge, and induced odd-cycle inequalities.
Moreover, every connected claw-free graph of even order has a perfect matching \cite{LasVergnas1975,Sumner1974}, and Plummer proved the following result.

\begin{theorem}[{\cite{Plummer1994}}]\label{t2}
Every $3$-connected claw-free graph of even order is bicritical.
\end{theorem}

The following theorem concerns the connectivity of $k$-factor-critical graphs.

\begin{theorem}[{\cite{Favaron1996,Yu1993}}]\label{t1}
Let $G$ be a $k$-factor-critical graph of order $n$.
Then $G$ is $k$-connected, $(k+1)$-edge-connected, and $(k-2)$-factor-critical whenever $k\geq 2$.
\end{theorem}

Theorem \ref{t1} implies that every $k$-factor-critical graph has minimum degree at least $k+1$.
In particular, every bicritical graph has minimum degree at least three.
For a graph $H$, let $c_o(H)$ denote the number of odd components of $H$.
A $k$-factor-critical graph $G$ is \emph{minimal} if $G-e$ is not $k$-factor-critical for every edge $e$ of $G$.
Moreover, for an edge $e=uv\in E(G)$,  $S_e\subseteq V(G)\setminus\{u,v\}$ is called a \emph{good set for $e$} if $|S_e|\geq k$, $c_o(G-S_e-e)=|S_e|-k+2$, and $u$ and $v$ belong to distinct odd components of $G-S_e-e$.
Favaron and Shi \cite{FavaronShi1998} gave the following characterization of minimal $k$-factor-critical graphs.

\begin{theorem}[{\cite{FavaronShi1998}}]\label{t4}
A $k$-factor-critical graph $G$ is minimal if and only if for each edge $e$ of $G$, there is a good set $S_e$.
\end{theorem}

A connected graph $G$ other than $K_2$ is \emph{matching covered} if every edge of $G$ is contained in a perfect matching.
An edge $e$ of a matching covered graph $G$ is \emph{removable} if $G-e$ is matching covered.
A matching covered graph $G$ is \emph{minimal} if $G$ has no removable edges.
Since every bicritical graph is matching covered, every bicritical graph that is minimal matching covered is also minimal bicritical.
Thus, the claw-free bicritical graphs that are minimal matching covered, characterized by Zhang et al. \cite{ZhangWangYuan2022}, are a subclass of claw-free minimal bicritical graphs.
Recently, Zhang et al. \cite{ZhangWangYuanNgCheng2025} gave a complete characterization of claw-free minimal matching covered graphs.
However, the classes of minimal bicritical graphs and minimal matching covered graphs are incomparable, even when restricted to claw-free graphs.
Examples are given in Section 2.
In this paper, we give complete structural characterizations of claw-free minimal factor-critical and claw-free minimal bicritical graphs.
Our main theorems are stated below.

\begin{theorem}\label{t3}
A graph $G$ is claw-free minimal factor-critical if and only if $G$ is an odd cycle or has at least two blocks, each of which is a triangle, with every cut vertex belonging to exactly two blocks.
\end{theorem}

Here, a \emph{block} of a graph $G$ is a maximal connected subgraph of $G$ having no cut vertex.
Clearly, a block is either a maximal $2$-connected subgraph or a $K_2$.
A cycle (resp. path) is called \emph{odd or even} according to whether its length is odd or even.
Let $C_n$ denote a cycle with $n$ vertices.
For $n\geq 3$, a \emph{wheel} $W_n$ is obtained from $C_n$ by adding a new vertex adjacent to every vertex of $C_n$.
We use $\overline{C_6}$ to denote the Cartesian product of $C_3$ and $K_2$, which is the complement graph of $C_6$.
For a graph $G$, replacing a vertex $v$ of degree three with a triangle consists of deleting $v$ from $G$, adding a new triangle, and distributing the three edges incident with $v$ to distinct vertices of the triangle.

To state our other result, we introduce the following operations.
Let $H$ be a simple $3$-connected cubic graph, and let $M$ be a matching of $H$ with $|M|$ even.
For each $v\in V(H)$, replace $v$ with a triangle $T_v$ with vertex set $\{v_e: e \text{ is an edge of } H \text{ incident with } v\}$.
For each edge $e=uv\in E(H)$, join $u_e$ to $v_e$ if $e\notin M$, and identify $u_e$ with $v_e$ if $e\in M$.
Denote the resulting graph by $\Theta(H,M)$.
Let $G_1$ and $G_2$ be two vertex-disjoint graphs other than $K_2$.
For each $i\in\{1,2\}$, let $x_iy_i\in E(G_i)$ satisfy
\textnormal{(i)} $d_{G_i}(x_i)=d_{G_i}(y_i)=3$,
\textnormal{(ii)} $G_i-x_iy_i$ is claw-free,
\textnormal{(iii)} $G_i-\{x_i,y_i\}$ is connected, and
\textnormal{(iv)} both $N_{G_i}(x_i)\setminus\{y_i\}$ and $N_{G_i}(y_i)\setminus\{x_i\}$ induce a $K_2$.
An edge satisfying conditions \textnormal{(i)}--\textnormal{(iv)} is called an \emph{admissible edge}.
Delete the admissible edges $x_1y_1$ and $x_2y_2$ from $G_1$ and $G_2$, respectively, and identify either $x_1$ with $x_2$ and $y_1$ with $y_2$, or $x_1$ with $y_2$ and $y_1$ with $x_2$.
The resulting graph, denoted $G_1\odot G_2$, is called an \emph{admissible $2$-sum} of $G_1$ and $G_2$.

\begin{theorem}\label{t7}
A graph $G$ is claw-free minimal bicritical if and only if it can be obtained from $K_4$, $W_5$, $\overline{C_6}$, and a graph $\Theta(H,M)$ by repeatedly applying admissible $2$-sums, where $H$ is a simple $3$-connected cubic graph, and $M$ is a matching of $H$ with an even number of edges.
\end{theorem}

The paper is organized as follows.
Section 2 introduces notation and preliminary results.
Sections 3 and 4 prove Theorems \ref{t3} and \ref{t7}, respectively.

\section{Preliminaries}
We use the following result to compare claw-free minimal bicritical graphs with claw-free minimal matching covered graphs.

\begin{lemma}[{\cite{Lovasz1987}}]\label{l2}
Every brick distinct from $K_4$ and $\overline{C_6}$ has a removable edge.
\end{lemma}

\begin{figure}[ht]
\centering
\begin{tikzpicture}[
    scale=0.45,
    line width=0.8pt,
    point/.style={
        circle,
        draw,
        fill=white,
        inner sep=1.6pt
    }
]

\coordinate (t1)  at (0,3.2);
\coordinate (t2)  at (2.6,3.2);
\coordinate (t3)  at (1.3,2.2);
\coordinate (t4)  at (1.29,1.05);

\coordinate (t5)  at (0.4,-0.55);
\coordinate (t6)  at (2.2,-0.55);

\coordinate (t7)  at (-2.1,-0.55);
\coordinate (t8)  at (-0.47,-1.12);
\coordinate (t9)  at (-0.63,-2.77);

\coordinate (t10) at (3.2,-1.12);
\coordinate (t11) at (4.7,-0.55);
\coordinate (t12) at (3.37,-2.77);

\draw
(t1)--(t2)
(t1)--(t3)
(t2)--(t3)
(t3)--(t4)
(t1)--(t7)
(t2)--(t11)
(t4)--(t5)
(t4)--(t6)
(t5)--(t6)
(t5)--(t8)
(t6)--(t10)
(t7)--(t8)
(t7)--(t9)
(t8)--(t9)
(t9)--(t12)
(t10)--(t11)
(t10)--(t12)
(t11)--(t12);

\node[point] at (t1) {};
\node[point] at (t2) {};
\node[point] at (t3) {};
\node[point] at (t4) {};
\node[point] at (t5) {};
\node[point] at (t6) {};
\node[point] at (t7) {};
\node[point] at (t8) {};
\node[point] at (t9) {};
\node[point] at (t10) {};
\node[point] at (t11) {};
\node[point] at (t12) {};

\end{tikzpicture}
\caption{The truncated tetrahedral graph.}
\label{fig1}
\end{figure}
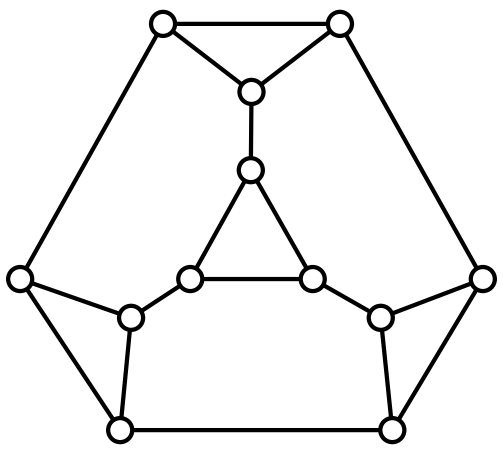

A cycle $C_6$ is claw-free minimal matching covered, but not bicritical, and so not minimal bicritical.
Let $T_{12}$ denote the truncated tetrahedral graph, which is obtained from $K_4$ by replacing each vertex with a triangle, see Figure \ref{fig1}.
Then $T_{12}$ is cubic, 3-connected and claw-free.
By Theorem \ref{t2}, $T_{12}$ is bicritical, and so is a brick.
For each edge $e=uv\in E(T_{12})$, both $u$ and $v$ have degree two in $T_{12}-e$.
So $T_{12}-e$ is not bicritical.
Therefore, $T_{12}$ is claw-free minimal bicritical.
Since $T_{12}$ is neither $K_4$ nor $\overline{C_6}$, by Lemma \ref{l2}, $T_{12}$ has a removable edge, and so is not minimal matching covered.
Thus we get the following proposition.

\begin{prop}
Neither the set of claw-free minimal bicritical graphs nor the set of claw-free minimal matching covered graphs contains the other.
\end{prop}

Let $H$ be a subgraph of $G$.
An \emph{ear} of $H$ in $G$ is an odd path in $G$ whose two distinct ends lie in $V(H)$ and whose internal vertices lie in $V(G)\setminus V(H)$.
An ear is \emph{trivial} if it is $K_2$.
An \emph{odd ear decomposition} of a graph $G$ is a sequence of subgraphs $(G_0,G_1,\ldots,G_r)$ such that $G_0=K_1$, $G_r=G$, and
$G_i=G_{i-1}+P_i$ for each $i\geq1$, where $P_1$ is an odd cycle and $P_i$ is an  ear of $G_{i-1}$ for $i\geq2$.

\begin{lemma}[{\cite{LovaszPlummer}}]\label{l3}
Let $G$ be a 2-connected factor-critical graph.
Then $G$ has an odd ear decomposition $(G_0,G_1,\ldots,G_r)$ such that $G_0=K_1$ and each $G_i$ is a 2-connected factor-critical subgraph of $G$ for $i\in\{1,\ldots,r\}$.
Moreover, if $G$ is minimal factor-critical, then each ear in an odd ear decomposition of $G$ is nontrivial.
\end{lemma}

\begin{lemma}[{\cite{LovaszPlummer}}]\label{l6}
A graph is minimal factor-critical if and only if it is connected and each of its blocks is minimal factor-critical.
\end{lemma}

\begin{lemma}[{\cite{LovaszPlummer}}]\label{l8}
A minimal bicritical graph contains a wheel as a subgraph if and only if it is a wheel.
\end{lemma}

A $2$-vertex cut $S$ of a graph $G$ is called a \emph{$2$-separation} if $c_o(G-S)\neq |S|$.
For a $2$-vertex cut $S$ of a bicritical graph $G$,  since $G-S$ has a perfect matching, we have $c_o(G-S)=0$, and so $S$ is a $2$-separation.

\begin{lemma}[{\cite{ZhangWangYuan2022}}]\label{l9}
Let $G$ be a bicritical graph with a $2$-separation $\{u,v\}$.
Then
\vspace{-6pt}
\begin{enumerate}[(i)]
\setlength{\topsep}{0pt}
\setlength{\partopsep}{0pt}
\setlength{\itemsep}{-1pt}
\setlength{\parsep}{0pt}
\item both $u$ and $v$ have at least two neighbours in each component of $G-\{u,v\}$;
\item if $uv\in E(G)$, then $G-uv$ is bicritical and $uv$ is removable in $G$.
\end{enumerate}
\end{lemma}

\begin{lemma} \label{11}
Let $G$ be the union of two graphs $G_1$ and $G_2$ that are connected graphs other than $K_2$ with $V(G_1)\cap V(G_2)=\{u,v\}$.
For $i=1,2$, let $H_i=G_i+uv$ if $uv\notin E(G_i)$, and $H_i=G_i$ otherwise.
Then the following hold.
\vspace{-6pt}
\begin{enumerate}[(i)]
\setlength{\itemsep}{-1pt}
\item \mbox{(\cite{LovaszPlummer})}\label{l10}
$G$ is bicritical if and only if both $H_1$ and $H_2$ are bicritical.
\item \mbox{(\cite{GuoWuZhang2024})}\label{l11}
If $G$ is bicritical and $\{u,v\}$ is a $2$-separation of $G$, then, for $i=1,2$,
\[
\mathcal{D}(G)\cap E(H_i)=
\begin{cases}
\mathcal{D}(H_i)\setminus\{uv\}, & \text{if } uv\notin E(G),\\[2pt]
\mathcal{D}(H_i)\cup\{uv\},     & \text{if } uv\in E(G).
\end{cases}
\]
\end{enumerate}
\end{lemma}

Here, $\mathcal{D}(G)=\{e\in E(G): G-e\text{ is bicritical}\}$.  Lemma \ref{11}(\ref{l10}) is contained in the proof of Lemma 5.2.8 in \cite{LovaszPlummer}.
The following lemma summarizes the structural properties
in Lemma 3.2 and Propositions 3.4 and 3.5 of
\cite{GuoLiLuZhang2024}.

\begin{lemma}\label{l12}
Let $G$ be a claw-free minimal k-factor-critical graph with $k\geq 2$, $e=uv\in E(G)$ with $d_G(u), d_G(v)\geq k+2$,  $S_e$ be a good set for $e$ of minimum cardinality, and $G_1,\ldots,G_t$ be the components of $G-e-S_e$ such that $u\in V(G_1)$ and $v\in V(G_2)$. Then $G-e-S_e$  has no even components and exactly one of the following holds:
\vspace{-6pt}
\begin{enumerate}[(i)]
\setlength{\itemsep}{-1pt}
\item $|S_e|=k$, $t=2$ and $|V(G_i)|\geq3$ for $i=1,2$.
(In this case, $e$ is called a \emph{type I edge}.)
\item $|S_e|=k+1$, $t=3$, and $V(G_1)=\{u\}$, $V(G_2)=\{v\}$. For each $w\in S_e$, $\{uw,vw\}\subseteq E(G)$ and $N_G(w)\cap V(G_3)\neq\varnothing$.
(In this case, $e$ is called a \emph{type II edge}.)
\end{enumerate}
\end{lemma}

By Lemma 2.8, we deduce the following corollary, which implies that there are no type II edges in the case $k=2$.

\begin{corollary}\label{cotype1}
Let $G$ be a claw-free minimal bicritical graph.
Then every edge $e=uv\in E(G)$ with
$d_G(u),d_G(v)\geq4$ is a \emph{type I edge}.
\end{corollary}

\begin{proof}
Since $G$ is minimal bicritical, by Theorem \ref{t4}, we have a good set for $e$.
Let $S_e$ be a good set for $e$ of minimum cardinality.
By Lemma \ref{l12} with $k=2$, the edge $e$ is either a type I edge or a type II edge.
Suppose that $e$ is a type II edge.
By Lemma \ref{l12}(ii), we have $|S_e|=3$, $t=3$, $V(G_1)=\{u\}$ and $V(G_2)=\{v\}$.
Moreover, every vertex of $S_e$ is adjacent to both $u$ and $v$.
Since $G$ is claw-free, $S_e$ is not an independent set.
Thus, there exist adjacent vertices $a,b\in S_e$.
Then $G[\{u,v,a,b\}]\cong K_4$.
By Lemma \ref{l8}, $G$ is a wheel.
However, a wheel has at most one vertex of degree at least four, contradicting $d_G(u),d_G(v)\geq4$.
Thus, $e$ cannot be a type II edge, and hence it is a type I edge.
\end{proof}

The following lemma follows from Claim 1 in the proof of Theorem 1.6 of \cite{GuoLiLuZhang2024} for $k=2$.

\begin{lemma}\label{l13}
Let $G$ be a $3$-connected claw-free minimal bicritical graph.
Then the subgraph of $G$ induced by the type I edges whose ends both have degree at least four is a forest.
\end{lemma}

\section{Proof of Theorem \ref{t3}}

\begin{lemma}\label{l7}
Every $2$-connected claw-free minimal factor-critical graph is an odd cycle.
\end{lemma}

\begin{proof}
Let $G$ be a $2$-connected claw-free minimal factor-critical graph.
By Lemma \ref{l3}, $G=P_0+P_1+\cdots+P_t$, where $P_0=K_1$, $P_1$ is an odd cycle, and $P_i$ is an ear of $G_{i-1}$ for each $i\geq2$, where $G_{i-1}=P_0+P_1+\cdots+P_{i-1}$.
We now show that $t=1$ by induction on the order of $G$.

Suppose, to the contrary, that $t\geq 2$.
By Lemma \ref{l3}, every $P_i$ with $i\geq2$ is nontrivial.
Let $P_t=x_0x_1\cdots x_{r}$.
Then $r\geq3$.
Since $G_{t-1}$ is an induced subgraph of $G$, it is claw-free.
By Lemma \ref{l3}, $G_{t-1}$ is $2$-connected and factor-critical.
If $G_{t-1}-e$ is factor-critical for some $e\in E(G_{t-1})$, then $G-e=(G_{t-1}-e)\cup P_t$ is factor-critical, contradicting the minimality of $G$.
Thus $G_{t-1}$ is minimal factor-critical.
By induction hypothesis, $G_{t-1}$ is an odd cycle.
So $t=2$.
Let $u$ and $v$ be the two neighbours of $x_0$ in $G_{1}$.
Since $G$ is claw-free, we have $uv\in E(G)$.
As $G_{1}$ is an induced subgraph of $G$, we have $uv\in E(G_{1})$, and so $G_{1}$ is a triangle $uvx_0u$.
Suppose, without loss of generality, that $x_{r}=u$.
Deleting the edge $x_0u$ from $G$ yields an odd cycle, which is factor-critical, again contradicting the minimality of $G$.
Therefore, $t=1$, and so $G$ is an odd cycle.
\end{proof}

\begin{proof}[Proof of Theorem \ref{t3}]
An odd cycle is clearly claw-free minimal factor-critical.
Suppose that $G$ has at least two blocks, every block is a triangle, and every cut vertex belongs to exactly two blocks.
Since $K_3$ is minimal factor-critical, Lemma \ref{l6} implies that $G$ is minimal factor-critical.
Since every vertex of $G$ lies in a triangle, $G$ is claw-free.
This proves the sufficiency.

Now, suppose that $G$ is a claw-free minimal factor-critical graph.
By Lemma \ref{l6}, every block of $G$ is minimal factor-critical.
Since $G$ is a nontrivial connected graph and $K_2$ is not factor-critical, every block of $G$ is $2$-connected.
Moreover, each block is an induced subgraph of $G$ and so is claw-free.
By Lemma \ref{l7}, each block of $G$ is an odd cycle.
If $G$ has only one block, then $G$ is an odd cycle.

Suppose that $G$ has at least two blocks, and let $B$ be a block of $G$.
Then $B$ contains a cut vertex $u$ that belongs to another block $B'$.
If $B$ has at least five vertices, then the two neighbours $u_1$ and $u_2$ of $u$ on $B$ are nonadjacent.
Let $u_3\in N_{B'}(u)$. Then $G[\{u,u_1,u_2,u_3\}]$ is a claw.
This contradiction implies that $B$ is a triangle.

If a cut vertex $v$ of $G$ belongs to at least three blocks, then choosing one neighbour of $v$ from each of three blocks, we get a claw of $G$.
Hence every cut vertex of $G$ belongs to exactly two blocks.
This completes the proof.
\end{proof}

\section{Proof of Theorem \ref{t7}}

For a graph $H$, let $V_4(H)=\{v\in V(H):d_H(v)=4\}$.
We write $V_4$ for $V_4(H)$ when no confusion can arise.
Let $\mathcal{G}$ denote the set of all $3$-connected claw-free minimal bicritical graphs.
We first characterize $\mathcal{G}$ and then extend the characterization to the general case.

\begin{lemma}\label{l15}
Let $G\in\mathcal{G}$ be a graph that is not a wheel.
Then $d_G(v)\in\{3,4\}$ for each $v\in V(G)$, $G[V_4]$ is a forest, and $G$ contains no $K_4$.
\end{lemma}

\begin{proof}
By Lemma \ref{l8}, $G$ contains no wheel and hence contains no $K_4$.
For each $v\in V(G)$, if there is a cycle in $G[N_G(v)]$, then this cycle together with $v$ forms a wheel.
So $G[N_G(v)]$ is a forest.
Since $G$ is claw-free, $G[N_G(v)]$ has no independent set of size three.
Since every forest on at least five vertices has an independent set of size at least three, we have $d_G(v)\leq4$.
Since $G$ is bicritical, we have $d_G(v)\geq3$.
So $d_G(v)\in\{3,4\}$.
Let $e=uv$ be an edge of $G[V_4]$.
Then $d_G(u)=d_G(v)=4$.
By Corollary \ref{cotype1}, $e$ is a type I edge.
By Lemma \ref{l13}, $G[V_4]$ is a forest.
\end{proof}

\subsection{Reduction to line graphs}

For a graph $H$, the \emph{line graph} $L(H)$ has vertex set $E(H)$, with two vertices adjacent if and only if the corresponding edges of $H$ have a common end.
For each $e=xy\in E(H)=V(L(H))$, the neighbour set of $e$ in $L(H)$ is the union of two cliques corresponding to the edges of $H$ incident with $x$ and $y$, respectively.
Hence, $L(H)$ is claw-free.

\begin{lemma}[{\cite{Harary}}]\label{l16}
Let $G$ be a claw-free graph.
Then $G$ is a line graph of a simple graph if and only if every pair of odd triangles of $G$ that share exactly one edge induces a $K_4$.
\end{lemma}

Here, a triangle $T$ in a graph $G$ is said to be \emph{odd} if there is a vertex in $V(G)\setminus V(T)$ that has one or three neighbours in $T$.
We use $\Delta(G)$ and $\delta(G)$ to denote the maximum degree and the minimum degree of $G$, respectively.
The \emph{square} of a graph $H$, denoted by $H^2$, is the graph with vertex set $V(H)$ in which two vertices are adjacent if and only if their distance in $H$ is at most two.
Recall that $C_n$ is a cycle with $n$ vertices.
The proof of the following lemma is inspired by an argument in \cite{BruhnStein2012}.

\begin{lemma}\label{l17}
Let $G$ be a $3$-connected claw-free $K_4$-free graph such that $\Delta(G)\leq4$ and $G[V_4]$ is a forest.
If $G$ is not a line graph of a simple graph, then $G\cong C_7^2-v$ for some $v\in V(C_7^2)$.
\end{lemma}

\begin{proof}
We first present the following observation, which follows from the fact that $G$ is claw-free and will be used repeatedly.

\begin{observation}
Let $u\in V(G)$, and $x, y, z, w\in N_G(u)$. If $xyz$ is an induced path in $G[N_G(u)]$, then either $xw\in E(G)$ or $zw\in E(G)$.
\end{observation}

We now show that $G$ contains the square of a path on six vertices as a subgraph.
Since $G$ is not a line graph, by Lemma \ref{l16}, there are two odd triangles $u_1u_2u_3$ and $u_2u_3u_4$ in $G$ sharing exactly the edge $u_2u_3$ and not inducing a $K_4$.
Hence $u_1u_4\notin E(G)$.
Since $G$ is $3$-connected, at least one of $u_2$ and $u_3$ has a neighbour not in $\{u_1,u_2,u_3,u_4\}$; otherwise, $\{u_1,u_4\}$ would be a $2$-vertex cut of $G$.
Let $u_5\in V(G)\setminus \{u_1,u_2,u_3,u_4\}$.
Suppose, without loss of generality, that $u_3u_5\in E(G)$.
Since $u_1u_2u_4$ is an induced path in $G[N_G(u_3)]$, by the above observation, we may assume, without loss of generality, that $u_4u_5\in E(G)$.
Since $G$ contains no $K_4$, we have $u_2u_5\notin E(G)$.
Since the triangle $T=u_2u_3u_4u_2$ is odd, and both $u_1$ and $u_5$ are adjacent to exactly two vertices of $T$, there is a vertex $u_6\in V(G)\setminus \{u_1,\ldots,u_5\}$ adjacent to an odd number of vertices of $T$.
Since $G$ contains no $K_4$, $u_6$ is adjacent to exactly one vertex of $T$.
Since $\Delta(G)\leq4$, we have $u_3u_6\notin E(G)$.
Suppose, without loss of generality, that $u_4u_6\in E(G)$.
Then $u_2u_6\notin E(G)$.
Since $u_2u_3u_5$ is an induced  path in $G[N_G(u_4)]$, the above observation implies $u_5u_6\in E(G)$.
Let $Q_6=u_1u_2u_3u_4u_5u_6$.
Then $Q_6^2$ is a subgraph of $G$.

Note that $d_G(u_3)=d_G(u_4)=4$ and $\Delta(G)\leq4$.
Thus, neither $u_3$ nor $u_4$ has a neighbor in $V(G)\setminus V(Q_6)$.
If either $u_2$ or $u_5$ has a neighbor in $V(G)\setminus V(Q_6)$, then either $\{u_2,u_3,u_4\}$ or $\{u_3,u_4,u_5\}$ consists of three vertices, each of degree $4$, that induce a triangle, contradicting that $G[V_4]$ is a forest.
If $V(G)\setminus V(Q_6)\neq\varnothing$, then $\{u_1,u_6\}$ is a $2$-vertex cut of $G$, contradicting that $G$ is $3$-connected.
So $V(G)=V(Q_6)$.

Recall that $Q_6^2$ is a subgraph of $G$. If $u_2u_5\in E(G)$, then $G[\{u_2,u_3,u_4,u_5\}]\cong K_4$, contradicting that $G$ is $K_4$-free.
Thus $u_2u_5\notin E(G)$.
Since $G$ is $3$-connected, $\delta(G)\geq 3$.
If $u_1u_6\notin E(G)$, then
$u_1u_5,u_2u_6\in E(G)$.
Thus $u_2,u_3$ and $u_4$  have degree $4$ and induce a triangle, contradicting that $G[V_4]$ is a forest.
Thus $u_1u_6\in E(G)$.
Consequently, up to isomorphism, $G$ is one of $C_6^2$, $C_6^2-u_1u_5$, or $C_6^2-u_1u_5-u_2u_6$.
In each of the first two cases, the vertices$u_2,u_3$ and $u_4$  have degree $4$ and induce a triangle, contradicting that $G[V_4]$ is a forest.
Hence $G\cong C_6^2-u_1u_5-u_2u_6\cong C_7^2-v_7$.
This completes the proof.
\end{proof}

\begin{lemma}\label{l18}
For each $v\in V(C_7^2)$, the graph $C_7^2-v$ is bicritical but not minimal bicritical.
\end{lemma}

\begin{proof}
Let $C_7=v_1v_2\cdots v_7v_1$.
Since $C_7^2$ is vertex-transitive, it suffices to consider $C_7^2-v_7$.
Since $C_7^2-v_7$ is $3$-connected and claw-free, by Theorem \ref{t2}, it is bicritical.
Since $(C_7^2-v_7)-v_3v_4\cong\overline{C_6}$ is bicritical, $C_7^2-v_7$ is not minimal bicritical.
\end{proof}

The only claw-free wheels of even order are $K_4$ and $W_5$.
By Lemmas \ref{l15}, \ref{l17}, and \ref{l18}, we obtain the following corollary.

\begin{corollary}\label{c1}
Every graph in $\mathcal{G}$ other than $K_4$ and $W_5$ is a line graph of a simple graph but is not a wheel.
\end{corollary}

\subsection{3-connected Graphs}

For a graph $H$ and a matching $M$ of $H$, let $\hat H(M)$ denote the graph obtained from $H$ by replacing each edge in $E(H)\setminus M$ with a path of length two.
For an edge $e\in E(H)$, let $P_e$ denote the corresponding edge or path of length two in $\hat H(M)$.

\begin{lemma}\label{l19}
Let $H$ be a simple $3$-connected cubic graph and $M$ a matching of $H$ with $|M|$ even.
Then the graph $\hat H(M)$ is simple and $\Theta(H,M)\cong L(\hat H(M))$.
\end{lemma}

\begin{proof}
For each edge $e=uv\in E(H)\setminus M$, $P_e$ is the path of length two in $\hat H(M)$ corresponding to $e$.
Let $x_e$ be the internal vertex of $P_e$.
Then $V(\hat H(M))=V(H)\cup\{x_e:e\in E(H)\setminus M\}$ and $E(\hat H(M))=M\cup\{ux_e,vx_e:e=uv\in E(H)\setminus M\}$.
Since $H$ is simple, so is $\hat H(M)$.

We define a map $\varphi:V(\Theta(H,M))\longrightarrow E(\hat H(M))$ as follows.
If $e=uv\notin M$, define $\varphi(u_e)=ux_e$ and $\varphi(v_e)=vx_e$.
If $e=uv\in M$, then $u_e$ and $v_e$ coincide in $\Theta(H,M)$.
Denote this vertex by $w$ and define $\varphi(w)=uv$.
Clearly, $\varphi$ is a bijection.
Note that every edge of $\Theta(H,M)$ that does not belong to a triangle corresponds to an edge in $E(H)\setminus M$.
For two adjacent vertices $u_e$ and $v_e$ with respect to an edge $e=uv\in E(H)\setminus M$, $\varphi(u_{e})$ and $\varphi(v_{e})$ are two edges of $\hat H(M)$ incident with $x_e$; for two vertices $v_e$ and $v_{e'}$ of $\Theta(H,M)$ belonging to the same triangle $T_v$,   $\varphi(v_{e})$ and $\varphi(v_{e'})$ are two edges of $\hat H(M)$ incident with $v$.
Therefore, the images under $\varphi$ of any two adjacent vertices in $\Theta(H,M)$ are adjacent in $L(\hat H(M))$.

Conversely, every pair of adjacent edges in $\hat H(M)$ shares an end that is either a vertex $v\in V(H)$ or an internal vertex $x_e$ of a path of length two.
Recall that $H$ is a cubic graph.
For the former case, the corresponding vertices of $\Theta(H,M)$ lie in the triangle $T_v$ and so are adjacent.
For the latter case, the two edges are $ux_e$ and $vx_e$ for some edge $e=uv\in E(H)\setminus M$, and their corresponding vertices $u_e$ and $v_e$ are adjacent in $\Theta(H,M)$.
Therefore, $\varphi$ is an isomorphism, and hence $\Theta(H,M)\cong L(\hat H(M))$.
\end{proof}

\begin{lemma}\label{l20}
For a graph $H$, if $L(H)$ is $3$-connected, then deleting at most two edges from $H$ leaves at most one component containing edges.
\end{lemma}

\begin{proof}
Let  $E'\subseteq E(H)$ with $|E'|\leq 2$.
Suppose, to the contrary, that $H-E'$ has two distinct components $H_1$ and $H_2$, each containing an edge.
Let $e_i\in E(H_i)$ for $i\in\{1,2\}$.
Since $L(H)$ is $3$-connected, deleting the vertices of $L(H)$ in $E'$ leaves a connected graph $L'$.
Hence, there is a path in $L'$ connecting the vertices $e_1$ and $e_2$.
This path yields a sequence of edges in $H-E'$, which lies in a common component of $H-E'$ and contains $e_1$ and $e_2$, contradicting that $e_1$ and $e_2$ lie in distinct components of $H-E'$.
\end{proof}

For a nonempty proper subset $X$ of $V(G)$, we denote by $\partial_G(X)$ the set of all the edges of $G$ with one end in $X$ and the other end not in $X$.

\begin{lemma}\label{l21}
Let $H$ be a graph such that $d_H(v)\in\{2,3\}$ for every $v\in V(H)$ and no two vertices of degree two are adjacent.
If $L(H)$ belongs to $\mathcal{G}$ but is not a wheel, then the edges of $H$ whose ends both have degree three form a matching.
\end{lemma}

\begin{proof}
Let $G=L(H)$.
Suppose, to the contrary, that there are two edges $e_1=uv$ and $e_2=vw$ in $H$, where $d_H(u)=d_H(v)=d_H(w)=3$.
Let $e_3=vz$ be another edge of $H$ incident with $v$.
Since $d_G(e_1)=d_G(e_2)=4$, by Corollary \ref{cotype1}, $e_1e_2$ is a type I edge.
By Lemma \ref{l12}(i), there exists $S\subseteq V(G)\setminus\{e_1,e_2\}$ with $|S|=2$ such that $G-e_1e_2-S$ has exactly two components, one containing $e_1$ and the other containing $e_2$.

Since $e_1e_3e_2$ is a path of length two in $G-e_1e_2$, we have $e_3\in S$.
Let $f$ be the other vertex in $S$, which is an edge of $H$.
Then $S=\{e_3,f\}$.
If $H-v-f$ contained a path $P$ joining $u$ and $w$, then $e_1$, the edges of $P$, and $e_2$ induce a path in $G-e_1e_2-S$ joining $e_1$ and $e_2$, contradicting that $e_1$ and $e_2$ lie in distinct components of $G-e_1e_2-S$.
Thus $u$ and $w$ lie in distinct components of $H-v-f$.

For the two components of $H-v-f$ containing $u$ and $w$, there is one that does not contain $z$, and we denote its vertex set by $X$.
Suppose, without loss of generality, that $u\in X$.
Then an edge in $\partial_H(X)$ is either $f$ or incident with $v$.
Since $w, z\notin X$, we have $\partial_H(X)\subseteq\{e_1,f\}$.
Since $d_H(u)=3$, at least one edge incident with $u$ lies in $H[X]$.
Note that $v, w \notin X$.
Thus, $H-\partial_H(X)$ has at least two components, one of which is $H[X]$ and has an edge, while another contains $e_2$, contradicting Lemma \ref{l20}.
The proof is complete.
\end{proof}

\begin{theorem}\label{t6}
A graph $G$ belongs to $\mathcal{G}$ if and only if $G$ is isomorphic to $K_4$, $W_5$, $\overline{C_6}$, or $\Theta(H,M)$, where $H$ is a simple $3$-connected cubic graph and $M$ is a matching of $H$ with $|M|$ even.
\end{theorem}

\begin{proof}
We first prove the sufficiency.
Since $K_4$, $W_5$, and $\overline{C_6}$ are $3$-connected claw-free graphs of even order, by Theorem \ref{t2}, they are bicritical.
Since every edge of each graph has an end of degree $3$, deleting any edge leaves a vertex of degree $2$.
Thus the resulting graph is not bicritical.
Consequently, all three graphs are minimal bicritical, and so belong to $\mathcal{G}$.

Let $H$ be a $3$-connected cubic graph and $M$ a matching of $H$ with $|M|$ even.
Let $G=\Theta(H,M)$ and $\hat H=\hat H(M)$.
By Lemma \ref{l19}, $G\cong L(\hat H)$, and so $G$ is claw-free.
Let $X\subseteq V(G)$ with $|X|\leq 2$, and let $E_X$ be the set of edges of $\hat H$ corresponding to $X$.
Then $G-X\cong L(\hat H-E_X)$.
Let $F_X=\{e\in E(H):E(P_e)\cap E_X\neq\varnothing\}$.
Then $|F_X|\leq 2$.
Since $H$ is $3$-edge-connected, $H-F_X$ is connected.
Replacing each edge $e$ of $H-F_X$ by the path $P_e$ in $\hat H$ yields a connected subgraph $H'$ of $\hat H-E_X$ containing all vertices of $H$.
Every edge of $\hat H-E_X$ not in $H'$ lies on $P_e$ for some $e\in F_X$ and has an end in $V(H)$.
Thus all edges of $\hat H-E_X$ induce a component.
Consequently, $L(\hat H-E_X)$ is connected, and so is $G-X$.
Therefore, $G$ is $3$-connected.

Since $H$ is cubic, $|V(H)|$ is even.
Recall that $|M|$ is even.
Thus $|V(G)|=3|V(H)|-|M|$ is even.
By Theorem \ref{t2}, $G$ is bicritical.
Note that every vertex of $G$ has degree three or four, and the vertices of degree four in $G$ are precisely those corresponding to the edges of $M$.
Since $M$ is a matching, the vertices of degree four in $G$ form an independent set of $G$.
Hence, deleting any edge of $G$ leaves a vertex of degree two, so the resulting graph is not bicritical.
Therefore, $G$ is minimal bicritical and hence belongs to $\mathcal{G}$.

We now prove the necessity.
Let $G\in\mathcal{G}$.
Suppose that $G$ is isomorphic to neither $ K_4$ nor $W_5$.
It suffices to prove that $G$ is isomorphic to $\overline{C_6}$ or $\Theta(H,M)$, where $H$ is a $3$-connected cubic graph and $M$ is a matching of $H$ with $|M|$ even.

By Corollary \ref{c1}, $G=L(H_1)$ for some simple graph $H_1$, and $G$ is not a wheel.
Since $G$ is 3-connected, we may choose $H_1$ to be connected.
By Lemma \ref{l15}, every vertex of $G$ has degree three or four, and $G$ contains no $K_4$.
If there are four edges incident with the same vertex of $H_1$, then they induce a $K_4$ in $L(H_1)$.
Thus $\Delta(H_1)\leq3$.
For each edge $xy\in E(H_1)$, we have $d_G(xy)=d_{H_1}(x)+d_{H_1}(y)-2\in\{3,4\}$.
It follows that $d_{H_1}(v)\in\{2,3\}$ for each $v\in V(H_1)$ and no two vertices of degree two in $H_1$ are adjacent.
Thus $H_1$ has some vertices of degree three.
By Lemma \ref{l21}, the edges of $H_1$ whose ends both have degree three form a matching of $H_1$, denoted by $M$.
For each vertex $x$ of degree two in $H_1$ with neighbors $u$ and $v$, deleting $x$ from $H_1$ and adding the edge $uv$, we get a graph $H$, which may have multiple edges and is connected, loopless, and cubic.
Clearly, $M\subseteq E(H)$, and $H_1=\hat{H}(M)$.
Let $n=|V(H)|$.
Since $H$ is cubic, $n$ is even.
Moreover, $|V(G)|=|E(H_1)|=|M|+2\bigl(|E(H)|-|M|\bigr)=3n-|M|$.
Since $G$ is bicritical, $|V(G)|$ is even, and hence $|M|$ is even.

We now show that $H$ is $3$-edge-connected.
Suppose that $H-F$ is disconnected for some $F\subseteq E(H)$ with $|F|\leq2$.
Since $H$ is loopless and cubic, every component of $H-F$ has an edge.
For each edge $e\in F$, choose an edge $f_e\in E(P_e)$, and let $F_{H_1}:=\{f_e: e\in F\}$.
Then $|F_{H_1}|\leq2$, and the subgraphs of $H_1$ corresponding to distinct components of $H-F$ lie in distinct components of $H_1-F_{H_1}$.
Hence $H_1-F_{H_1}$ has at least two components containing edges.
Recall that $L(H_1)$ is 3-connected.
We get a contradiction to Lemma \ref{l20}.
Thus $H$ is $3$-edge-connected.

If $H$ is simple, then $H$ is $3$-connected since $H$ is cubic and $3$-edge-connected. By Lemma \ref{l19}, $G=L(H_1)= L(\hat H(M))\cong\Theta(H,M)$.
If $H$ is not simple, let $u$ and $v$ be two vertices of $H$ that are joined by parallel edges.
Since $H$ is 3-edge-connected and cubic, we have $V(H)=\{u,v\}$.
Thus $H$ consists of three parallel edges joining $u$ and $v$.
Since $M$ is an even matching, we have $M=\varnothing$, and hence $H_1\cong K_{2,3}$ and $G=L(H_1)\cong\overline{C_6}$.
This completes the proof.
\end{proof}

\subsection{Admissible $2$-sums}

\begin{lemma}\label{l22}
An admissible $2$-sum of two claw-free minimal bicritical graphs is claw-free minimal bicritical.
\end{lemma}

\begin{proof}
For $i=1, 2$, let $G_i$ be a claw-free minimal bicritical graph, and let $e_i=x_iy_i$ be an admissible edge of $G_i$.
Set $H_i=G_i-e_i$ for $i=1,2$. Then $H_i$ is claw-free, connected and not isomorphic to $K_2$.
Let $G=G_1\odot G_2$ be the admissible $2$-sum of $G_1$ and $G_2$ with respect to $e_1$ and $e_2$.
Then $E(G)=E(H_1)\cup E(H_2)$.
Suppose, without loss of generality, that $G$ is obtained from $H_1$ and $H_2$ by identifying $x_1$ and $x_2$ as a vertex $x$ and $y_1$ and $y_2$ as a vertex $y$.

Since $G_1$ and $G_2$ are bicritical, by Lemma \ref{11}(\ref{l10}), $G$ is bicritical.
Since each $G_i$ is minimal bicritical, we have $\mathcal D(G_i)=\varnothing$.
By Lemma \ref{11}(\ref{l11}), we have $\mathcal D(G)=\varnothing$, and hence $G$ is minimal bicritical.

Suppose that $G$ contains an induced claw with center $v$.
Since both $H_1$ and $H_2$ are claw-free, we have $v\in\{x,y\}$.
Note that the neighbors of each of $x$ and $y$ in $G$ induce a graph consisting of two disjoint copies of $K_2$.
Thus neither $x$ nor $y$ can be the center of an induced claw, a contradiction.
Therefore, $G$ is claw-free.
\end{proof}

\begin{lemma}\label{l23}
Let $G$ be a claw-free minimal bicritical graph that is not $3$-connected, and let $\{x,y\}$ be a $2$-separation of $G$.
Then:
\vspace{-6pt}
\begin{enumerate}[(i)]
\setlength{\topsep}{0pt}
\setlength{\partopsep}{0pt}
\setlength{\itemsep}{-1pt}
\setlength{\parsep}{0pt}
\item $xy\notin E(G)$;
\item $G-\{x,y\}$ has exactly two components, each of even order;
\item $d_G(x)=d_G(y)=4$, and each of $x$ and $y$ has exactly two neighbours in each component of $G-\{x,y\}$;
\item if $J$ is a component of $G-\{x,y\}$, then $G[N_G(x)\cap V(J)]\cong K_2$, $G[N_G(y)\cap V(J)]\cong K_2$.
\end{enumerate}
\end{lemma}

\begin{proof}
By Lemma \ref{l9}\textup{(ii)} and the minimality of $G$, we have $xy\notin E(G)$.
\textup{(i)} follows.

Since $\{x,y\}$ is a $2$-separation of $G$, each component of $G-\{x,y\}$ has even order.
By Lemma \ref{l9}\textup{(i)}, each of $x$ and $y$ has at least two neighbours in each component of $G-\{x,y\}$.
Suppose, to the contrary, that $G-\{x,y\}$ has at least three components.
Then $x$, together with one of its neighbours from each of three distinct components, induces a claw of $G$.
This contradiction proves \textup{(ii)}.

Since every wheel is $3$-connected, $G$ is not a wheel.
By Lemma \ref{l8}, $G$ contains no wheel.
Since $G$ is bicritical, we have $\delta(G)\geq3$.
Hence $G[N_G(v)]$ is a forest for each $v\in V(G)$.
Since $G$ is claw-free, $G[N_G(v)]$ has no independent set of size $3$.
Every forest on at least five vertices has an independent set of size at least $3$, so $\Delta(G)\leq4$.
Since $xy\notin E(G)$ and each of $x$ and $y$ has at least two neighbours in each of the two components of $G-\{x,y\}$, we have $d_G(x)=d_G(y)=4$, and each has exactly two neighbours in each component.
This proves \textup{(iii)}.

By (ii), $G-\{x,y\}$ has two components.
Then any two neighbours of $x$ in the component $J$ of $G-\{x,y\}$ are adjacent since $G$ is claw-free.
Recall that $G[N_G(v)]$ is a forest.
Combining  (iii), we have $G[N_G(x)\cap V(J)]\cong K_2$.
The same argument applies to $y$, proving (iv).
\end{proof}

\begin{lemma}\label{l24}
Let $G$ be a claw-free minimal bicritical graph that is not $3$-connected.
Then $G$ is an admissible $2$-sum of two claw-free minimal bicritical graphs.
\end{lemma}

\begin{proof}
By Theorem \ref{t1}, $G$ is $2$-connected.
Since $G$ is bicritical and not $3$-connected, it has a $2$-separation, say $\{x,y\}$.
By Lemma \ref{l23}\textup{(ii)}, $G-\{x,y\}$ has exactly two components, say $J_1$ and $J_2$, each of even order.
For $i=1,2$, let $G_i=G[V(J_i)\cup\{x,y\}]$ and $H_i=G_i+xy$.
Then $G_i$ is connected and not isomorphic to $K_2$.
By Lemma \ref{l23}\textup{(i)}, $xy\notin E(G)$, and so $H_i$ is simple.
By Lemma \ref{11}(\ref{l10}), $H_i$ is bicritical.

Since $G_i$ is an induced subgraph of $G$, it is claw-free.
So, any induced claw of $H_i$ must use the edge $xy$, and hence its center is $x$ or $y$.
However, the two neighbours of each of $x$ and $y$ in $J_i$ are adjacent by Lemma \ref{l23}\textup{(iv)}.
Thus $H_i$ is claw-free, and $H_i-xy=G_i$ is claw-free.

We next show that $H_i$ is minimal bicritical.
By Lemma \ref{l23}\textup{(iii)}, $d_{G}(x)=d_{G}(y)=4$, and so $d_{H_i}(x)=d_{H_i}(y)=3$.
Then $H_i-xy$ has a vertex $x$ of degree two.
This implies that $H_i-xy$ is not bicritical.
Therefore, $xy\notin\mathcal D(H_i)$.
Since $G$ is minimal bicritical, $\mathcal D(G)=\varnothing$.
By Lemma \ref{11}(\ref{l11}), no edge of $G_i$ belongs to $\mathcal D(H_i)$.
Hence $\mathcal D(H_i)=\varnothing$ as required.

It is easy to see that $xy$ is an admissible edge of $H_i$ for $i=1,2$.
Hence $G$ is an admissible $2$-sum of the two claw-free minimal bicritical graphs $H_1$ and $H_2$.
\end{proof}

\begin{proof}[Proof of Theorem \ref{t7}]
Let $\mathcal B$ be the set of graphs obtained from $K_4$, $W_5$, $\overline{C_6}$, and a graph $\Theta(H,M)$ by repeatedly applying admissible $2$-sums, where $H$ is a simple 3-connected cubic graph, and $M$ is a matching of $H$ with an even number of edges.
By Theorem \ref{t6}, $\mathcal G\subseteq \mathcal B$.
Recall that $\mathcal{G}$ is the set of all $3$-connected claw-free minimal bicritical graphs.
By Lemma \ref{l22}, every graph in $\mathcal B$ is claw-free minimal bicritical.

Conversely, let $G$ be a claw-free minimal bicritical graph.
We prove that $G\in\mathcal B$ by induction on $|V(G)|$.
If $G$ is 3-connected, then $G\in\mathcal G\subseteq\mathcal B$.
Suppose that $G$ is not $3$-connected.
By Lemma \ref{l24}, $G$ is an admissible $2$-sum of two claw-free minimal bicritical graphs $G_1$ and $G_2$.
By the induction hypothesis, $G_1,G_2\in\mathcal B$, and so $G\in\mathcal B$.
\end{proof}

\section*{Acknowledgements}
This work is supported by the National Natural Science Foundation of China (Nos. 12571381, 12501499
and 12371361) and Natural Science Foundation of Henan Province (No. 262300421874).

\section*{Declaration of AI Use}
AI was used during the exploratory stage of this paper. All mathematical arguments and proofs in the
manuscript were developed and rigorously verified by the authors.

\end{document}